\documentclass[11pt]{article}

\usepackage[latin1]{inputenc}
\usepackage{amsfonts,amsmath,amssymb,amsthm,graphicx,epsfig,float}
\usepackage[T1]{fontenc}

\usepackage[english,francais]{babel}

{
  
  \newtheorem*{theoreme*}{Th\'eor\`eme}

\newtheorem*{corollaire*}{Corollaire}
\newtheorem*{proposition*}{Proposition}
\theoremstyle{remark}
  \newtheorem*{remarque*}{Remarque}
}

\newcounter{ex}

\newenvironment{rem*}{
  \noindent\textbf{Remarque. }}{}

\newcommand{\Cc}{\mathbb{C}}

\renewcommand {\epsilon}{\varepsilon}

\renewcommand {\leq}{\leqslant}
\renewcommand {\geq}{\geqslant}

\title{{\bf Une démonstration du théorème de recouvrement de surfaces d'Ahlfors}}
\author{Henry de Thélin}
\date{}

\begin{document}
\maketitle

%$\Cc, \Rr$%

%% Redefinition Titre
\def\figurename{{Fig.}}%
\def\proofname{Preuve}% for AMS-\LaTeX
\def\contentsname{Sommaire}%
%% Fin

\begin{abstract}

Nous donnons une preuve courte du théorème de recouvrement de surfaces d'Ahlfors.
\end{abstract}

\selectlanguage{english}
\begin{center}
{\bf{A proof of Ahlfors' theorem on covering surfaces }}
\end{center}

\begin{abstract}

We give a short proof of Ahlfors' theorem on covering surfaces.

\end{abstract}

\selectlanguage{francais}

Mots-clefs: topologie, surfaces, applications holomorphes.\\
AMS: 30D35, 30C25.

\section*{{\bf Introduction}}
\par

Soit $f: \Sigma \longrightarrow \Sigma_0$ une application holomorphe
non constante entre surfaces de Riemann compactes, connexes et qui ont
éventuellement du bord.

Si $f$ est propre, de degré $d$, la formule de Riemann-Hurwitz donne
la caractéristique d'Euler de $\Sigma$ en fonction de celle de
$\Sigma_0$. Plus précisément:
$$\chi(\Sigma) + r = d \chi(\Sigma_0),$$
où $r$ est le nombre de ramifications de $f$. En particulier, nous avons l'inégalité:
$$\chi(\Sigma) \leq d \chi(\Sigma_0).$$
Le théorème d'Ahlfors (voir \cite{A}) permet d'étendre en quelque sorte cette inégalité au cas non propre. Plus précisément, si $S$ désigne le nombre de feuillets moyen au-dessus de $\Sigma_0$, i.e.: $$S=\frac{\mbox{aire de }f(\Sigma) \mbox{ comptée avec multiplicité }}{\mbox{aire de }\Sigma_0},$$
et $L$, la longueur du bord relatif de $\Sigma$ (i.e. la longueur de
$f(\partial \Sigma) - \partial \Sigma_0$):

\begin{theoreme*}{\label{Ahlfors}}
On a l'inégalité:
$$\min (\chi(\Sigma),0) \leq S \chi(\Sigma_0)+hL,$$ 
où $h$ est une constante qui ne dépend que de $\Sigma_0$.

\end{theoreme*}

L'objectif de cet article est de donner une démonstration courte de ce
théorème qui a permis à L. Ahlfors de réaliser une version géométrique
de la théorie des distributions de valeurs de R. Nevanlinna.

Signalons aussi que ce théorème a servi récemment dans des problèmes de dynamique holomorphe (voir par exemple \cite{BLS}).

\section*{{\bf Démonstration du théorème}}
\par

Commençons par préciser dans quel cadre nous allons nous placer.

Le théorème d'Ahlfors n'a d'intérêt que si $\chi(\Sigma_0)$ est
négatif (chose que l'on supposera dans la suite). Ensuite, dans ce
théorème, la surface $\Sigma_0$ peut avoir du bord et on peut
considérer des métriques très générales (voir par exemple
\cite{N}). Cependant, bien que notre démonstration s'étende à ces cas,
par souci de simplicité, nous supposerons d'une part que $\Sigma_0$
n'a pas de bord et d'autre part qu'elle est munie d'une métrique
lisse.

Fixons maintenant un point $p$ de $\Sigma_0$. En utilisant la
décomposition canonique de $\Sigma_0$ (voir \cite{R} p.23), nous
pouvons construire $2g$ courbes simples et lisses (où $g$ est le genre
de $\Sigma_0$) qui passent par le point $p$, qui ne s'intersectent pas
en dehors du point $p$, et telles que la surface $\Sigma_0$ privée de
ces courbes soit un disque $D_0$. Ce disque possède $4g$ côtés que
nous pouvons supposer ne passant pas par les valeurs critiques de
$f$.

Notons $\mathcal{C}$, l'ensemble des composantes connexes de $\Sigma$
au-dessus du disque $D_0$ (i.e. dans $f^{-1}(D_0)$).

Voici maintenant le plan de la suite de la démonstration. Dans un
premier paragraphe, nous allons enlever les composantes $\Delta$ de $\mathcal{C}$ pour
lesquelles $f(\Delta)$ est très différent de $\Sigma_0$. Cette
modification de $\Sigma$ permettra de se rapprocher du cas où $f$ est
un revêtement. Le second paragraphe sera ensuite consacré au calcul de
$\chi(\Sigma)$.

\subsection*{{\bf 1) Modification de $\Sigma$}}

\par

Si $\Delta$ est une composante au-dessus de $D_0$, on peut décomposer
$f(\Delta)$ en strates: $\Delta^{1}$ est la partie de $f(\Delta)$
recouverte au moins une fois,..., $\Delta^{p}$ celle recouverte $p$
fois (où $p$ est la multiplicité maximale de $f(\Delta)$). On note
$\gamma_j$ la partie de $f(\partial \Sigma)$ qui borde $\Delta^{j}$ et
on fixe $\epsilon > 0$ petit (en particulier devant la longueur du
plus petit côté de $D_0$). Alors,
si $l(\gamma_j)$ désigne la longueur de $\gamma_j$, on peut avoir:

- l'existence d'un $j$ dans $\{1,...,p \}$ avec $l(\gamma_j) \geq
  \epsilon$ (on dira que la strate en question a un bord long)

ou:

- pour tout $j$ dans $\{1,...,p \}$, $l(\gamma_j) \leq \epsilon$ (on
  parlera de bord court).

En utilisant une inégalité isopérimétrique linéaire et le fait que le
bord de $D_0$ est lisse par morceaux, on remarque que dans le dernier
cas, on a:\\
Soit:
$$ \mbox{aire de } \Delta^j \geq \mbox{aire de } D_0 - h l(\gamma_j) \mbox{   , pour un certain } j \in \{1,...,p \}.$$
Soit:
$$\mbox{aire de } \Delta^j \leq h  l(\gamma_j) \mbox{   , pour tout }
j \in \{1,...,p \},$$
où $h$ est une constante qui ne dépend que de $\Sigma_0$. Dans toute
la suite ce type de constante sera toujours notée $h$.

Dans notre contexte, ce sont les composantes du dernier type
$(\Sigma^i)_{i \in I}$ qui ont leur image la plus différente de
$\Sigma_0$. Les $\Sigma^i$ seront appelées mauvaises composantes et
nous allons les enlever.

Remarquons que nous avons deux possibilités:

soit $\Sigma= \cup_{i \in I}\Sigma^i$, ce qui implique:
$$S \leq hL$$
par définition même des mauvaises composantes (dans ce cas le théorème
est démontré),

soit $\Sigma \neq \cup_{i \in I}\Sigma^i$ et alors:
$$\chi(\Sigma) \leq \chi(\Sigma - \cup_{i \in I}\Sigma^i),$$
car les composantes $\Sigma^i$ (avec $i$ dans $I$) ont du bord en
commun avec $\Sigma$. 

Par ailleurs, toujours par définition des mauvaises composantes, la
surface $\widetilde{\Sigma}=\Sigma - \cup_{i \in I}\Sigma^i$ a la
longueur de son bord relatif majoré par $hL$ et l'aire de
$f(\widetilde{\Sigma})$ diffère de celle de $f(\Sigma)$ d'au plus
$hL$.

En résumé, si on démontre que 
$$\chi(\widetilde{\Sigma}) \leq \widetilde{S} \chi(\Sigma_0)+h \widetilde{L}$$
(avec des notations évidentes), le théorème sera démontré. Dans toute
la suite nous oublierons les tildes: $\Sigma$ désignera la surface modifiée.
 
\subsection*{{\bf 2) Calcul de $\chi(\Sigma)$}}

\par

Désignons toujours par $\mathcal{C}$ l'ensemble des composantes connexes de $\Sigma$
au-dessus du disque $D_0$ (i.e. dans $f^{-1}(D_0)$). On notera $s$ le nombre
de relevés du point $p$ par $f$. Parmi les composantes connexes dans
les préimages par $f$ des côtés de $D_0$ qui se trouvent dans la
surface $\Sigma$, il y en a certaines qui ne sont pas dans le bord de $\Sigma$. Nous les appellerons arêtes de $\Sigma$ et nous
noterons $a$ leur nombre. La surface $\Sigma$ privée des $s$ relevés
du point $p$ s'obtient en recollant les composantes de $\mathcal{C}$
le long des $a$ arêtes. En particulier, nous avons:
$$ \chi(\Sigma) = \sum_{\Delta \in \mathcal{C}}{\chi(\Delta)}-a+s.$$

Grâce à notre simplification effectuée au paragraphe précédent, on
sait qu'une composante de $\Sigma$ au-dessus du disque $D_0$ possède
une strate qui a soit un bord long, soit une aire supérieure à
$\mbox{aire}(D_0) - hl(\gamma_j)=A_0-hl(\gamma_j)$.

Il y en a au plus $\frac{L}{\epsilon}$ qui possèdent une strate
avec un bord long.

Pour les composantes qui possèdent une strate d'aire supérieure à
 $A_0-hl(\gamma_j)$, on voit facilement qu'il y en a au plus
 $S+hL$. On obtient alors la majoration:
$$\sum_{\Delta \in \mathcal{C}}{\chi(\Delta)} \leq S+hL.$$
Ensuite, si on choisit le point $p$ de sorte à être peu recouvert, nous avons $s \leq S$.
Pour finir, il nous faut donc minorer le nombre d'arêtes $a$.

On se fixe une composante connexe $\Delta$ au-dessus de $D_0$.

La composante $f(\Delta)$ contient des strates pour lesquelles on a
$l(\gamma_j) \geq \epsilon$ et d'autres avec cette longueur plus
petite que $\epsilon$. Les strates de la deuxième espèce se divisent
en deux catégories: celles qui ont une aire inférieure à
$hl(\gamma_j)$ et celles pour qui cette aire est supérieure à
$A_0-hl(\gamma_j)$. Si $m(\Delta)$ désigne le nombre de strates de la
dernière sorte, on a:

$$\nu(\Delta) \geq 4g m(\Delta),$$ 
où $\nu(\Delta)$ désigne le nombre d'arêtes de $\Sigma$ qui bordent
$\Delta$. En effet, soit $\gamma$ un sous-segment d'un côté de $D_0$
choisi de sorte qu'un de ses petits voisinages dans $D_0$ soit inclus
dans les $m(\Delta)$ strates (c'est possible car elles s'emboîtent les
unes dans les autres et que $\epsilon$ est très petit devant la
longueur du plus petit côté de $D_0$);
alors $\gamma$ a au moins $m(\Delta)$ relevés par $f$ qui sont inclus
dans des arêtes de $\Delta$.

La minoration de $\sum m(\Delta)$ entraîne donc celle du nombre
d'arêtes $a \geq \frac{1}{2} \sum \nu(\Delta)$. Pour l'obtenir, on va
minorer l'aire recouverte par les $\sum m(\Delta)$ strates.

Les strates pour lesquelles la longueur de $\gamma_j$ est supérieure à
$\epsilon$ sont en nombre au plus égal à $\frac{L}{\epsilon}$. L'union
de ces éléments est donc d'aire inférieure à
$\frac{L}{\epsilon}\mbox{aire}(D_0)=hL$.

De même, l'ensemble des strates qui ont un petit bord et une aire
majorée par $hl(\gamma_j)$ a une aire majorée par $hL$.

En combinant ces deux majorations, on obtient que l'aire recouverte par les strates qui ont un petit bord et une aire minorée par  $A_0-hl(\gamma_j)$ est supérieure à:
$$A_0S-hL.$$
Autrement dit, on a:
$$a \geq \frac{4g}{2} \sum_{\Delta \in \mathcal{C}} m(\Delta) \geq \frac{4g}{2} \frac{1}{A_0} (A_0S-hL) \geq 2gS-hL.$$
Finalement, on a:
$$\chi(\Sigma)= \sum_{\Delta \in \mathcal{C}}{\chi(\Delta)}-a+s \leq (1-2g+1)S+hL= \chi(\Sigma_0)S +hL,$$
qui est l'inégalité que l'on voulait démontrer.

\bigskip

Henry de Thélin\\
Université Paris-Sud (Paris 11)\\
Mathématique, Bât. 425\\
91405 Orsay\\
France

\end{document}